\input epsf
\title{Berwald metrics constructed by Chevalley's polynomials}

\newcommand{\cmt}[1]{\ifhmode\newline\fi{\sf *** \ \ #1 \\}}

\documentclass[11pt]{article}
\usepackage{latexsym}
\usepackage{epsfig}

\newtheorem{theorem}{Theorem}[section]
\newtheorem{lemma}[theorem]{Lemma}

\newtheorem{prop}[theorem]{Proposition}
\newtheorem{corol}[theorem]{Corollary}

\def\:{\colon}

\long\def\onefigure#1#2{
\begin{figure*}[tbh]
\begin{center}
#1
\end{center}
\caption{#2}
\end{figure*}
} 

\newcommand{\iipefig}[1]  
{\smallskip\begin{center}\def\IPEfile{#1.ipe}\input{\IPEfile}\end{center}\smallskip}
\newcommand{\lipefig}[2]  
{\onefigure{\def\IPEfile{frh-#1.ipe}\input{\IPEfile}}{\label{f:#1} #2} }

\author{
{\sc Z. I. Szab\'o}\thanks{Partially supported by
NSF grant DMS-0604861}
}  
\date{}

\begin{document}

\maketitle
\begin{abstract}

\noindent Berwald metrics are particular Finsler metrics
which still have linear Berwald connections. Their complete
classification is
established in an earlier work, \cite{sz1}, of this author. The
main tools in these classification are the 
Simons-Berger holonomy theorem
and the Weyl-group theory.
It turnes out that any Berwald metric is a perturbed-Cartesian product
of Riemannian, Minkowski, and such non-Riemannian metrics which can be
constructed on irreducible symmetric manifolds of $rank > 1$. 
The existence
of these metrics are well established by the above theories.

\noindent The present paper has several new
features. First, the Finsler functions of Berwald manifolds
are explicitly described by the
Chevalley polynomials. New results
are also the complete lists of reversible ($d(x,y)=d(y,x)$) resp.
irreversible ($d(x,y)\not =d(y,x)$) Berwald metrics. The Cartan
symmetric Finsler manifolds are also completely determined.
The paper is concluded by proving that
a Berwald metric is uniquely determined by the Minkowski metric induced
on an arbitrarily fixed maximal totalgeodesic 
flat submanifold (Cartan flat).
Moreover, any two Cartan flats are isometric.

\end{abstract}

\section{Introduction}

Loosly speaking, Finsler metrics generalize the  
Riemannian metrics like the Banach spaces do the Hilbert spaces.
On Riemannian manifolds the norm of tangent vectors is defined by 
a field, $\langle\, ,\, \rangle_{/p}$, of inner products, while
on a Finsler manifold, $(M,L)$,
the norm is defined by an
appropriate Banach-Minkowski norm 
$||X||_{/p}=L_p(X)$ on each tangent space $T_p(M)$. The
Finsler functions, $L$, correspond to the Lagrange functions
in the more general Lagrange-Hamilton theory. 
Analogously to the
Riemannian connections,
Berwald introduced the
so called Berwald connections 
on Finsler manifolds, which are,
however, non-linear in general. This
implies that the parallel
transports $\tau :T_p(M)\to T_q(M)$ 
defined along the curves of $M$ are just 
non-linear homogeneous maps preserving the norm $L_p(X)$.

A Finsler metric is called Berwald metric if the Berwald connection
is still a linear connection.
This concept is the closest possible to the concept 
of Riemannian metrics.
The prototype of Berwald metrics are the Minkowski
metrics defined by an appropriate norm on the affine space 
$\mathbf{R}^n$. 
In this case
the Berwald connection is nothing but the natural flat 
linear connection on the affine space.

The classification of Berwald metrics is established in 
an earlier paper (\cite{sz1}) of this author.
The basic tools used in this classification 
are the Simons-Berger holonomy theorem (asserting that the 
holonomy group $\mathcal{H}_p$
at a point $p\in M$ of an irreducible Riemannian manifold is transitive
on the unit sphere $I_p$ (indicatrix), unless the metric is 
symmetric of $rank\geq 2$), the Weyl group theory on symmetric spaces, 
and the deRham decomposition theory. 

It turned out that the linear Berwald
connection $\nabla$ on a Berwald manifold is always Riemann metrizable,
furthermore, the holonomy group, $\mathcal{H}_p$, 
leaves the indicatrix defined
by the Finsler function $L_p(X)$ invariant. Therefore, it remained 
the question that which Riemannian connections are metrizable also by 
non-Riemannian Finsler metrics. This metrization problem can be decided
at an arbitrarily fixed point $p\in M$ by means of the holonomy group
$\mathcal{H}_p$, since it boils down to the determination of the
Riemannian holonomy
groups which leave also non-Euclidean 
Finsler norms $||X||_{/p}=L_p(X)$ invariant.
Then the pursued Berwald metrics are nothing but the parallel extensions
of these norms onto the whole manifold $M$.
 
From the Simons-Berger theorem we get that for  
irreducible linear connections only  
the Riemannian connections of symmetric spaces with $rank\geq 2$
can be Berwald
connections of non-Riemannian Berwald metrics. 

The existence of the
desired Berwald metrics on symmetric manifolds is established
by showing that the above $\mathcal{H}_p$-invariant norms are uniquely
determined by their restrictions $L_{/\mathit a}$ onto an
arbitrarily fixed Cartan
subalgebra $\mathit a\subset T_p(M)$, furthermore, 
these restricted norms
are nothing but the Finsler norms which 
are invariant under the action
of the Weyl group $W_{\mathit a}$ acting on 
$\mathit a$. 
More precisely, exactly
these Weyl group invariant norms 
$L_{/\mathit a}$ have unique extensions, first,
to $\mathcal{H}_p$-invariant norms $L_p$ on the tangent space
$T_p(M)$, and then, 
by parallel displacement, to the desired Berwald metrics defined on
the whole manifold $M$. 
Thus, the question of existence is traced
back to the existence of the Weyl group invariant Finsler norms on
an arbitrarily fixed Cartan subalgebra.
Since the Weyl groups are finite, there is an
infinite dimensional variety of
the metrics appropriate to the problem.

The general reducible cases can be treated by the
deRham decomposition theorem. Let it also be mentioned that for
non-Riemannian Berwald metrics this decomposition means 
Cartesian product regarding
the manifold $M$ and the connection
$\Delta$, however, it is not the usual Cartesian 
product with respect to the metric. 
The usual Cartesian product of Finsler metrics
leads to singular (non-strictly convex) metrics in general.
However, one can easily overcome these difficulties by certain
perturbation of the metric. This technique is used both in \cite{sz1}
and this paper.  

Besides the new results, also 
the results accomplished in \cite{sz1} are reestablished 
in this paper in a
novel form. The key concepts used in these considerations 
are the 
{\it rank}, the {\it Cartan subspace}, the {\it Cartan flat},
and the {Weyl group} acting on a
Cartan subspace. These concepts will be 
introduced for general Riemannian
connections.

The {\it rank} of an irreducible locally symmetric Riemannian metric
(or a Riemann metrizable linear connection) is defined in the usual
way and the rank of an irreducible locally non-symmetric Riemannian
metric is defined by 1. For reducible Riemannian metrics the rank
is defined by the sum of ranks defined by the irreducible factors in
the de Rahm decomposition. 
Though it strongly relates to them,
this rank-concept is different from
those introduced in \cite{b,bbe,bbs}.

With respect to the holonomy group $\mathcal{H}_p$, 
consider the irreducible decomposition  
\[
T_p(M)=V_0\oplus S_1\oplus\dots\oplus S_l\oplus U_1\oplus\dots\oplus U_k
\]
of the tangent space at $p$, 
where $V_0$ is the maximal subspace fixed point-wise, 
and the subspaces $S_i$ resp. 
$U_i$ belong to irreducible symmetric resp. 
non-symmetric spaces (cf. (\ref{irrdec})).
Then a {\it Cartan flat},
$\mathit c$, is defined by the direct sum
\[
\mathit c=V_0\oplus \mathit{c}_{s1}\oplus\dots\oplus\mathit{c}_{sl}
\oplus
\mathit{c}_{u1}\oplus\dots\oplus\mathit{c}_{uk},
\]
where $\mathit{c}_{si}$ is a Cartan subalgebra in $S_i$ and 
$\mathit{c}_{ui}$ is an
arbitrary 1-dimensio\-nal subspace in $U_i$ (cf. (\ref{Cart})). 

The Weyl group $W_0$
on $V_0$ consists only of the identity map and 
the Weyl group $W_{si}$ acting
on the Cartan subalgebra $\mathit{c}_{si}$ of a 
symmetric Lie algebra is defined
by the standard definition. The Weyl group
$W_{ui}$ acting on $\mathit{c}_{ui}$ contains 
only two elements, the identity map
and the halfturn about the origin. Then the {\it Weyl group} 
$W_{\mathit c}$ acting on
$\mathit{c}$ is defined by the direct product (cf. (\ref{W})):
\[
W_{\mathit c}=W_0\times W_{s1}\times
\dots\times W_{sl}\times
W _{u1}\times\dots\times W_{uk}.
\]

The main results of this paper are
as follows. The basic question about characterizing the
non-Riemannian Berwald metrization is answered by
\begin{theorem}
{\bf (Abstract Main Theorem 1)} 
A torsion-free linear connection, $\nabla$, 
is Berwald
met\-rizable by a non-Riemannian Berwald metric if and only if
the closure, 
$\overline{\mathcal{H}}_p$, 
of the holonomy group is
compact (i. e., it is Riemann metrizable) 
and the rank of $\nabla$ is greater than 1.
\end{theorem}

There are two important corollaries to this theorem: 
\begin{theorem}
{\bf (Abstract Main Theorem 2)} 
(1) In the irreducible case, exactly the
Riemannian connections of symmetric spaces of $rank\geq 2$ are
Berwald metrizable by non-Riemann\-ian Berwald metrics. 

(2) Every
reducible Riemannian connection is
metrizable by non-Rie\-mann\-ian Berwald metrics. 
\end{theorem}

These statements can be reformulated as follows

\begin{theorem}
{\bf (Abstract Structure Theorem)} 
Let $(M,\nabla ,L)$ be a connected, simply connected, 
and complete Berwald
manifold. Then the affine manifold $(M,\nabla )$ decomposes into
the Cartesian product
$(M=M_0\times M_1\times\dots\times M_k ,
\nabla=\nabla_0\times \nabla_1\times\dots\times \nabla_k)$,
where 
$(M_0,\nabla_0,L_0)$
is the maximal Minkowskian factor and the irreducible factors 
$(M_i,\nabla_i,L_i)$
are either Riemannian manifolds or non-Riemannian affine symmetric
Berwald manifolds. The metric $L$ is a holonomy-invariant
(perturbed Cartesian) product
of the metrics $L_j$ defined on the factor manifolds.
\end{theorem}

By these theorems, the non-Riamannian factors
in the deRham decomposition of a Berwald manifold are either Minkowski
spaces or irreducible affine symmetric Berwald spaces of $rank\geq 2$.
Therefore, the Cartan classification of Riemannian symmetric spaces 
provides classification opportunity also for classifying Berwald 
manifolds. The complete list of irreducible
simply connected globally symmetric spaces can be found in Table 3.2.
One should exclude only the rank-one symmetric spaces (listed in 
(\ref{rank1})
in order to get the complete list of simply connected complete 
irreducible Riemannian spaces whose Riemannian connection can be
metrized also by non-Riemannian Berwald metric. The groups involved are
the unit components of the isometry groups also for the Berwald metrics.

One of the new features in this paper
is the explicit constructions of Berwald metrics.
As we have seen, all these metrics can be constructed by picking a 
$W_{\mathit c}$-invariant Finsler norm $L_{\mathit c}$ on an arbitrarily
fixed Cartan subspace $\mathit c$. The theory of polynomials invariant
under the action of Weyl groups (finite reflection groups) 
was developed by Chevalley
\cite{ch}. According to this theory, 
this polynomial ring is finitely generated
and all these invariant 
polynomials extend to unique
$\mathcal{H}_p$-invariant polynomials on the tangent space $T_p(M)$.

The explicit
construction of the Berwald metrics is 
achieved by means of these Chevalley 
polynomials. This approach describes the Berwald metrics in a
much more precise form. In the earlier paper
the investigations focused mainly on proving the existence of 
non-Riemannian Berwald metrics. In
 {\bf Constructive Main Theorem 1;2}
both the symmetric (reversible) and non-symmetric
(irreversible) B-metrics are completely classified.
 
These classifications are followed by two new topics.
The first one is the construction of all {\it 
Cartan-symmetric Finsler manifolds},
defined by those manifolds where 
there exists an
ivolutive isometry $\sigma_p$, 
for any point $p\in M$, 
such that $p$ is an isolated fix point 
of $\sigma_p$.
Then we prove:
 \begin{theorem}
A Finsler manifold is Cartan-symmetric if and only if it is a Berwald
manifold with a symmetric Berwald connection and absolute 
homogeneous metric.
\end{theorem}
This theorem along with
the classification of Berwald spaces provides the classification
of Cartan-symmetric Finsler manifolds.

In the other new topic the maximal totally geodesic
flats (called Cartan flats) 
on Berwald manifolds are investigated. These flats are defined 
by the surfaces inscribed by the 
geodesics tangent to a Cartan subspace $\mathit c$. 
We show that each vector is tangent to 
a Cartan flat and
the induced Minkowski metrics are
locally isometric on these flats.
 
By the above construction we get 
that the whole Berwald metric is uniquely
determined by the metric on an arbitrarily fixed
Cartan flat. 

Then we prove that
the Berwald connection of a Berwald space is
symmetric if and only if the isometries from the unit component $J^0$
of isometries act transitively on the set of pairs 
$(p,\mathit t)$, where $\mathit t$ is
a Cartan flat and $p\in \mathit t$. One can read this statement as 
follows: The k-flat homogeneous Berwald spaces are nothing but 
the symmetric Berwald spaces of rank $k$.

Note that this theorem generalizes the 
characterization of 
irreducible 
2-point homogeneous spaces
as being exactly the  
rank-one-symmet\-ric spaces.  
 
Berger's holonomy theorem is one of the central issues today, both in 
mathematics and mathematical physics. 
This theorem along with the Simons-Berger transitivity theorem was
introduced to me, 
at the time when the first paper \cite{sz1} was written, 
by J\'anos Szenthe. It is my pleasure to thank also \'Ad\'am
Kor\'anyi for the recent helpful conversations.

\section{Technicalities}

\noindent{\bf Berwald connections.}
A {\it Finsler manifold} is denoted by $(M,L)$, where $M$ is
the underlying manifold and $L$ is the smooth 
{\it Finsler function} defined on
the punctured tangent bundle $T(M)\backslash 0$. 
A coordinate neighborhood 
$x=(x^1,\dots ,x^n)$ on $M$ defines a coordinate neighborhood
$(x=(x^i),y=(y^i))$ on $T(M)$ such that $y^i=\dot{x}^i$ holds. Then,
with respect to the variable $y$, the 
$L(x,y)$ is a
{\it positive homogeneous function 
of degree 1 (called in short $+1$-homogeneous)},
furthermore, the
fundamental tensor 
$g_{ij}=(1/2)\dot{\partial}_i\dot{\partial}_jL^2$ is positive
definite. This latter assumption 
guarantees that the indicatrix $I_p$,
defined by vectors of unit norm in $T_p(M)$, is strongly convex.

Finsler functions having the symmetricity condition $L(x,y)=L(x,-y)$ are
called {\it absolute homogeneous}. They define symmetric metrics 
satisfying $d(x,y)=d(y,x)$. The more general metrics defined by
{\it positive homogeneous Finsler functions}
are non-symmetric, $d(x,y)\not =d(y,x)$, in general. 
A geodesic $\sigma (t), 0\leq t\leq r$, on a 
Finsler manifold is called reversible if
the reversed curve $\sigma (r-t)$ is also a geodesic. 
The non-reversible geodesics are
called also {\it irreversible ones}. It is well known 
that all the geodesics are reversible if the metric is 
defined by an absolute homogeneous Finsler function,
or, it is a Berwald metric \cite{bcs}. In the following
reversible resp. irreversible metrics are the synonyms for the
symmetric resp. non-symmetric ones defined by absolute resp.
proper positive homogeneous Finsler functions.
 
The {\it connection coefficients} 
$G^i,G^i_j,G^i_{jk}$ {\it of the Berwald 
connection} are defined by
\[
G^i=\frac{1}{4}g^{ik}(y^r\dot{\partial}_k\partial_rL^2-\partial_kL^2)
\,\, ,\,\,
G^i_j=\dot{\partial}_jG^i\,\, ,\,
\, G^i_{jk}=\dot{\partial}_jG^i_k=G^i_{kj} ,
\]
which are derived from the Euler-Lagrange equation 
$\ddot{x}^i+2G^i(x(t),\dot{x}(t))=0$ 
for a geodesic $x(t)=(x^i(t))$.
The coefficients $G^i_{jk}(x,y)$ 
correspond to the Christoffel symbols
on Riemannian manifolds. 
On general Finsler manifolds these symbols
are $+0$-homogeneous functions, which non-trivially depend on the
$y$-variable in general.

The {\it differential equation for 
a parallel vector field} $(y^i(t))$ 
along a curve $(x^i(t))$ is 
\[
{dy^i\over dt}+G^i_k(x(t),y(t)){dx^k\over dt}=0.
\]
This {\it parallel transport} keeps 
the Finsler-norm of the vectors, however,
it is non-linear if the functions $G^i_{jk}$ depend also on the 
$y$-variables. If $\nabla_k$ stands for the covariant derivative
with respect to the Berwald connection, then 
the {\it norm-keeping property}
is equivalent to the fundamental equation 
$\nabla_kL=L_{|k}=0$ of Finsler geometry.

One can give also an {\it invariant 
description of Berwald connections by means of
non-linear connections}. A general non-linear connection is 
defined by a smooth field (distribution), 
$H(x,y)$, of n-dimensional horizontal subspaces on the tangent bundle
$T(M)$, such that the subspace $H(x,y)$ 
is complement to the vertical subspace
at each point $(x,y)\in T(M)\backslash 0$. The connection coefficients, 
$G^i_j(x,y)$, on a coordinate system $(x^i,y^i)$ are defined such that
the field $E_j=\partial_j-G^i_j\dot{\partial}_i$ is horizontal. 
Among these general
non-linear connections the
Berwald connection of a Finsler manifold $(M,L)$ can be pinned down by
the following three properties:
\begin{prop} 
On any Finsler manifold there exists a uniquely
determined non-linear, so called Berwald connection, 
which is characterized by
the following properties: 

(1) The functions $G^i_j(x,y)$ are 
$+1$-homogeneous, 

(2) it is torsion free 
(i. e., $G^i_{jk}=G^i_{kj}$), and 

(3) the parallel transport defined by it is 
norm preserving (i. e., $\nabla_j L=E_j(L)=0$).
\end{prop}
This theorem is well known 
(cf. one of its demonstrations in\cite{sz1}).
\medskip

\noindent{\bf Basic tool applied in the constructions.}
A Finsler metric is said to be {\it Berwaldian} 
if the Berwald connection
is a linear connection, i. e., the connection coefficients $G^i_{jk}$
depend only on the $x$-variable. Let $(M,L,\nabla)$ 
be a Berwald manifold. Assume that $M$ is both arc-wise and simply
connected. Then, 
at an arbitrary point $p\in M$,
the holonomy group $\mathcal{H}_p$ of the linear connection
$\nabla$ 
has only one component.

By the above Proposition, the $\mathcal{H}_p$ 
is a subgroup of the compact
group, $G$, of linear transformations {\it leaving the indicatrix $I_p$
invariant.} Let $dg$ be the normalized invariant Haar measure 
on $G$ and consider 
an arbitrary inner product $\langle\, , \,\rangle^*$ 
on the tangent space $T_p(M)$.
Then the averaged inner product 
\begin{equation}
\langle X,Y\rangle =\int_G\langle g(X),g(Y)\rangle^*dg
\end{equation}
is invariant under the action both of $G$ and $\mathcal{H}_p$.
Therefore, by the parallel transports defined by $\nabla$, 
it can be extended into a Riemannian metric, $g(X,Y)$, for which 
the Riemannian
connection is $\nabla$. Thus {\it the linear connection $\nabla$ of a 
Berwald metric is Riemann metrizable} \cite{sz1}. 
By the above considerations we 
also have: {\it A torsion-free linear connection is Riemann metrizable
if and only if the closure $\overline{\mathcal H}_p$ of the holonomy 
group is compact.}

Thus the actual question is this: 
{\it Which Riemann connections can be
metrized also by non-Riemannian Berwald metrics?} This 
question {\it can be decided at an arbitrary fixed point $p$} by means
of the holonomy group. A subtler reformulation of the problem 
is this: {\it Which Riemannian holonomy groups $\mathcal{H}_p$ 
leave also non-Euclidean Finsler norm $||X||_p=L_p(X)$ invariant? 
After finding these norms, the parallel
extension of $L_p$ onto the whole manifold provides 
the desired Berwald metric.} 

The {\it rank} of an irreducible locally symmetric Riemannian metric
is defined in the usual
way and the rank of an irreducible locally non-symmetric Riemannian
metric is defined by 1. For reducible Riemannian metrics the rank
is defined by the sum of ranks of the irreducible factors appearing in
the deRahm decomposition. One of the main results of this paper is
the following
\begin{theorem}{\bf (Abstract Main Theorem 1.)} 
A torsion-free linear connection, $\nabla$, is Berwald
met\-ri\-zable by a non-Riemannian Berwald metric if and only if
the closure, 
$\overline{\mathcal{H}}_p$, 
of the holonomy group is
compact (i. e., the $\nabla$ is Riemannian metrizable)
and the rank of $\nabla$ is greater than 1.
\end{theorem}

There are two important corollaries to this theorem:
\begin{corol} {\bf (Abstract Main Theorem 2.)} 
(1) In the irreducible case, exactly the
Riemannian connections of symmetric spaces of $rank\geq 2$ are
Berwald metrizable by non-Riemannian Berwald metrics. 

(2) Every
reducible Riemannian connection is
metrizable by non-Rie\-mann\-ian Berwald metrics. 
\end{corol}

The name of the above statements indicates that they will be 
established in a much stronger form, namely, by explicit constructions
described in the next sections.

\section{Constructing the irreducible Berwald metrics}
{\bf Basic constructions.}
By the Berger-Simons theorem, the only candidates 
for irreducible Riemann
connections which can be metrized also 
by non-Riemannian Berwald metrics
are the irreducible symmetric connections of $rank\geq 2$.
Next, besides pointing out the existence of these Berwald metrics 
a technique is developed by which all these metrics can 
be constructed.  

In \cite{sz1}, the construction was established by 
Lemmas 2, 3, and 4. The
proofs are carried out by considering the {\it symmetric Lie algebra}
$T_p\oplus h_p$, at an arbitrarily fixed point 
$p\in M$ (cf. \cite{kn}, chapter 11). 
This decomposition of the Lie algebra corresponds to the
usual Cartan decomposition $\mathit g=\mathit p\oplus\mathit k$, 
i. e., the horizontal subspace is identified with
the tangent space $T_p(M)$ and the vertical subspace is identified with
the Lie algebra, $h_p$, of the holonomy group $exp(h_p)=\mathcal{H}_p$.
(On irreducible symmetric spaces the group $\mathcal H_p$ is 
identified with the unit
component, $K$, of isometries. Thus also the identification 
$h_p=\mathit k$ is well defined.) 

The Lie bracket $[T_p,T_p]\to h_p$ is defined by $[X,Y]=-R_p(X,Y)$,
where $R_p(X,Y)Z$ denotes the Riemannian curvature at $p$. On $h_p$
we keep the original Lie bracket, furthermore, $[h_p,T_p]=h_p(T_p)$.

Let $\mathit a$ be a maximal Abelian subalgebra 
({\it Cartan subalgebra})
in $T_p$. 
On geometric level they correspond to the tangent spaces of the 
{\it maximal totally geodesic flats.} 
It is well known that any two Cartan
subalgebras can be transformed to each other  
by isometries represented in $\mathcal{H}_p$.
Furthermore, every vector $v\in T_p$ 
is contained in at least one of the
Cartan subspaces. More precisely, there exist 
an everywhere dense open
subset of the so called {\it regular vectors} which belong
exactly to one of the Cartan subspaces. 
The others (called singular vectors) are covered by more Cartan
subspaces. The rank is defined by the common 
dimension of the Cartan subspaces.

Denotation $\mathcal{H}_p(X)$ stands for the 
orbit of the holonomy group
passing through $X\in T_p$. Then the tangent 
space to this orbit at $X$
is $h_p(X)$. The following statements are standard whose proofs 
can be found also in \cite{sz1}.
\begin{lemma} (A) Any Cartan subalgebra $\mathit a$ is 
intersected by any above
orbit at finite many points: 
$\mathit a\cap\mathcal{H}_p(X)=\{q_1,\dots ,q_l\}\not =\emptyset$.

(B) For every pair $q_i,q_j$ of intersection points there
exist some transformations $g\in \mathcal{H}_p$ 
leaving $\mathit a$ invariant
($g(\mathit a)=\mathit a$) and satisfying $g(q_i)=q_j$.

(C) Let $\mathcal W$ be the subgroup of $\mathcal{H}_p$, consisting of 
transformations leaving $\mathit a$ invariant. The $W$ stands for the
restriction of $\mathcal{W}$ onto $\mathit a$. Then $W$ is a finite
irreducible group acting transitively on the intersection
points $\{q_1,\dots ,q_l\}$, with respect to any orbit.
\end{lemma}

The $W$ is nothing but the 
{\it Weyl group} defined on a
Cartan subspace $\mathit a$. 

According to the construction technique
described in the previous section, we need to construct all the 
$\mathcal{H}_p$-invariant Finsler norms, 
$L_p(X)$, defined
on the tangent space 
$T_p(M)$ at a fixed point $p\in M$. The above statements
clearly reveal the following

\begin{corol}{\bf (Basic Observation)}
Any $\mathcal{H}_p$-invariant F-norm, 
$L_p(X)$, is uniquely determined
by its restriction, $L_{/\mathit a}(V)$, 
to an arbitrary Cartan subalgebra
$\mathit a\subset T_p(M)$. More precisely, the restricted 
(Finsler!) norms $L_{/\mathit a}$
are invariant under the action of the Weyl group $W$, furthermore,
$L_p(X)=L_{/\mathit a}(V_X)$ holds, 
where $V_X$ is an arbitrary vector from the 
intersection
$\mathit a\cap \mathcal{H}_p(X)$. 
\end{corol}

The question remained after this Observation is that which
$W$-invariant Fins\-ler norms (defined only on $\mathit a$) 
extend to
$\mathcal{H}_p$-invariant Finsler norms (defined on $T_p(M)$). 
Keep in mind that the extended norm should be both smooth
on the punctured space $T_p(M)\backslash 0$ and strongly
convex. The most simple technique
providing all the desired $\mathcal{H}_p$-invariant Finsler norms 
is the averaging of Finsler norms, $L^*_p(X)$, by the
holonomy group $\mathcal{H}_p$. This averaged norm,
\begin{equation}
L_p(X)=\int_{\mathcal{H}_p}L_p^*(g(X))dg,
\end{equation}
is obviously a Finsler norm. This construction does not use the
above Basic Observation and it gives only vague idea about the size
of the class of
$\mathcal{H}_p$-invariant Finsler norms. In order 
to discover the wide range of
the solutions 
(in fact, we prove
that the non-equivalent solutions form an
infinite dimensional variety) 
we turn to the technique offered by
Observation 3.2. Namely, we try to reach the solutions by extending
suitable $W$-invariant norms on an arbitrarily fixed Cartan 
subalgebra $\mathit a\subset T_p(M)$ to the whole tangent space
$T_p(M)$. This subtler construction technique exploits the theory
of Chevalley's polynomials.
\medskip 

\noindent{\bf Constructing by the Chevalley extension.}
If one considers an arbitrary 
$W$-invariant Finsler norm, $L_{/\mathit a}$, then only 
the continuity and convexity
of the extended norm $L_p(X)$ can be established easily. 
Surprisingly enough, the smoothness
and strong convexity lead to difficult problems.
These hardships are clearly indicated in the following review
of Weyl group invariant
functions on a subalgebra $\mathit a$.

The $W$-invariant polynomials on a Cartan subspace $\mathit a$ were
investigated by Chevalley \cite{ch}. In one of his theorems 
he proved:
{\it If a finite group $W$ of linear transformations 
of an $n$-dimensional
vector space is generated by reflections (the 
Weyl groups fall into this category!) 
then the algebra of polynomials invariant
under the action of $W$ is generated by $n$ algebraically 
independent homogeneous
polynomials and the identity.} In an other theorem he proved that
{\it any $W$-invariant polynomial $P$ on 
$\mathit a$ can be extended to
$Ad(K)=Ad(\mathcal{H}_p)$-invariant polynomial 
$\tilde{P}$ on $\mathit p=T_p(M)$}
(unpublished, cf. \cite{he}, p. 430). 

By the invariant theory developed
in \cite{gl}, \cite{sch}, \cite{ma}, {\it 
every $C^{\infty}$ $W$-invariant
function, $f$, on $\mathit a$ can be written as}
\[
f(x)=g(\sigma_1(x),\dots ,\sigma_n(x)),\quad g\in 
C^{\infty}(\mathbf{R}^n),
\]
{\it where $\sigma_1,\dots ,\sigma_n$ 
generate the $W$-invariant polynomials,
therefore, the $f$ extends to the $C^\infty$, $\mathcal{H}_p$-invariant
function $\tilde{f}$ on $T_p(M)$.} However, 
this statement is false if $C^\infty$
is replaced by $C^k$. This does not terminates, of course, the
possibility for the corresponding $C^k$-extension, however, 
as far the author knows, this
$C^k$-extension is still an open problem. 
Let it be also mentioned that a direct
elementary proof of the $C^\infty$-extension theorem can be found
in \cite{da}. First we describe a simple technique using this
$C^\infty$-extension.

The whole set of $W$-invariant $C^\infty$ Finsler norms on 
$\mathit a$ can be 
constructed by the averaging 
$L_{/\mathit a}(H)=\sum_{\varphi\in W}
L^*_{/\mathit a}(\varphi (H))/|W|$
of arbitrary $C^\infty$ Finsler norms $L^*_{/\mathit a}$ defined 
on $\mathit a$. 
Let $L_{/\mathit a}(H)$ be such an invariant norm. 
The $C^\infty$ Chevalley extension can not be directly
applied to this function, since it is of class $C^\infty$ only on the
punctured space $\mathit a\backslash 0$. The function 
$e^{-1/|H|^2}
L_{/\mathit a}(H)$, however,
is $C^\infty$ everywhere and its Chevalley extension is of the
form 
$e^{-1/|X|^2}
\tilde{L}_{/\mathit a}(X)$,
where
$\tilde{L}_{/\mathit a}(X)$ 
is the Chevalley extension of
$L_{/\mathit a}(H)$ 
onto $T_p(M)$. 

The next problem arising about the extended function is
that it is not strongly convex in general. One can control 
this problem by
applying an appropriate correction (perturbation) term 
by considering the norm
\begin{equation}
L(X)=(\gamma |X|^2+\tilde{L}^2(X))^{1/2},
\label{L1}
\end{equation}
where $|X|$ is a $\mathcal H_p$-invariant Euclidean norm and 
the constant $\gamma >0$ is chosen such that $L$ should be strongly
convex. Since the matrix field $\tilde{g}_{ij}=(1/2)\partial_i
\partial_j\tilde L^2$ is
$+0$-homogeneous, the existence
$\gamma$ is obvious. In fact, there exists a $\delta >0$ such that all
the constants satisfying $\gamma > \delta$ are 
appropriate to the problem. One can construct all the $C^\infty$
solutions by this technique.

\noindent{\bf Remark.}
A deeper insight reveals that any
positive $\gamma$ is suitable to the problem since the
matrix field is positive semi-definite everywhere. (The proof of
this latter statement is omitted).
\medskip

\noindent{\bf Chevalley polynomials.} In this section we briefly review
the theory describing Chevalley's polynomials explicitly. 
We are particularly interested in 
cases when there
are also skew-symmetric polynomials among the $W$-invariant
polynomials. One can construct non-symmetric (irreversible) 
Berwald metrics exactly in these cases. 

The Chevalley polynomials are described with the help of the root
systems of symmetric Lie algebras. Bourbaki \cite{bou} defines
an {\it abstract root system}, $\Delta$, in a real inner product space
$V$ by a finite set of non-zero elements of $V$ satisfying
the following properties:

(i) $\Delta$ spans $V$,

(ii) the orthogonal transformations
$
s_{\alpha}(\varphi )=\varphi 
-{2\langle\varphi ,\alpha\rangle\over |\alpha |^2}\alpha
$,
for all $\alpha\in\Delta$, 

leave $\Delta$ stable, and

(iii) $2\langle\beta ,\alpha\rangle/|\alpha |^2$ 
is an integer number for all
$\alpha ,\beta\in\Delta$.

\noindent It is well known that 
$-\alpha\in\Delta$ if 
$\alpha\in\Delta$. Furthermore, the only possible members of $\Delta$
proportional to $\alpha$ are $\pm\alpha$, or 
$
\pm\alpha$ and $\pm 2\alpha
$,
or
$
\pm\alpha$ and $\pm {1\over 2}\alpha
$.

An abstract root system is called {\it reduced} if $\alpha\in \Delta$
implies that $2\alpha$ is not in $\Delta$. The irreducible reduced
root systems are classified by means of the Dynkin diagram \cite{bou},
\cite{wa}.
According to this classification, the different classes are
\begin{equation}
A_l,B_l,C_l,(BC)_l,D_l,E_6,E_7,E_8,F_4,G_2,
\label{class}
\end{equation}
where the indexes are always equal to $dim(V)$.

The restricted root system, $\Delta_{\mathit p}$, on the dual space 
$V=\mathit{a}^*$ of a Cartan subalgebra $\mathit a\subset\mathit p$
of a Riemannian symmetric Lie algebra 
$\mathit g=\mathit p\oplus\mathit k$ is defined by the real non-zero
functionals $\alpha (H)$ for which there exist algebra elements
of the form $E_{\alpha}=X_{\alpha}+Z_{\alpha}$, where 
$X_{\alpha}\in\mathit p$ and $Z_{\alpha}\in\mathit k$, such that 
$[H,X_{\alpha}]=\alpha (H)Z_{\alpha}$,
$[H,Z_{\alpha}]=\alpha (H)X_{\alpha}$, and therefore
$[H,E_{\alpha}]=\alpha (H)E_{\alpha}$ hold.
Vectors $E_{\alpha}$ span the  
subspaces $\mathit{g}_\alpha$
in the {\it root-space decomposition} 
$\mathit{g}^{\mathbf C}=\mathit a^{\mathbf C}\oplus
\sum_{\alpha\in\Delta}\mathit{g}_{\alpha}$ resp.
$\mathit{g}=\mathit a\oplus
\sum_{\alpha\in\Delta}\mathit{g}_{\alpha}$ corresponding to the 
compact resp. non-compact cases. 
These restricted root systems are  always reduced
and, in the irreducible cases, each of them 
falls into one of the categories
described in (\ref{class}).

The Weyl group is generated by the reflections $s_\alpha$, where
$\alpha\in\Delta_{\mathit p}$, in the hyperplanes defined by 
$\alpha (H)=0$ in $\mathit a$. 
Therefore, it is determined by the root system. It
is well known \cite{bou}, \cite{wa} that the systems $B_l,C_l,(BC)_l$
determine the same Weyl group $\mathbf B_l$ and the
complete list of distinct irreducible Weyl groups is
\begin{equation}
\mathbf A_l,\mathbf B_l,\mathbf D_l,\mathbf E_6,\mathbf E_7,
\mathbf E_8,\mathbf F_4,\mathbf G_2.
\label{Wgroups}
\end{equation}

The Weyl group is transitive on a complete set
of roots having the same length. 
Thus there exists an isometry $w_{\alpha}\in W$
sending $\alpha$ to $-\alpha$, for all $\alpha\in\Delta$. 
This property is not satisfied for the other elements of 
$\mathit a^*$ in general. More precisely, it is satisfied for all vector
$v\in \mathit a^*$ if and only if the function space invariant under the
action of the Weyl group consists only even (reversible) functions.
This case can be characterized also as follows.
\begin{prop}
For any $v\in\mathit a$ 
there exist $w_v\in W$ 
reversing $v$ to $-v$, i. e.
$w_v(v)=-v$, 
if and only if 
$-id\in W$. 

There exist non-trivial skew (odd) functions, 
satisfying $\varphi (-H)=-\varphi (H)$, among the $W$-invariant
functions if and only if $W$ excludes the central symmetry $-id$.
\end{prop}

This statement can be 
established by considering the Weyl chambers, which are the 
connected components of the regular vectors $v$ of $\mathit a$ 
(recall that $v$ is regular if $\alpha (v)\not =0$ holds for all 
$\alpha\in\Delta$). By a standard statement of the theory,
the $W$ is simply transitive on the Weyl chambers. Fix
a Weyl chamber and let
$w\in W$ be the group element transforming this Weyl
chamber to the opposite one. If 
$-id\not\in W$ 
then there exist 
vectors $v$ in the chamber such that $w(v)\not =-v$ hold. 
By the 
simply transitivity of the 
Weyl group on the Weyl chambers, relation
$\varphi (v)\not =-v$ must hold 
for all $\varphi\in W$. This proves the existence of 
irreversible $v$'s, if $-id\not\in W$. 
If $-id\in W$, then each $v\in \mathit a$ is
obviously reversible. This proves the statement completely. 

The fundamental weights 
$\{\lambda_1,\dots ,\lambda_l\}$
with respect to a basis
$\{\alpha_1,\dots ,\alpha_l\}$ 
of the root system are defined by 
$2\langle\lambda_i,\alpha_j\rangle 
/\langle\alpha_j,\alpha_j\rangle=\delta_{ij}$ on $\mathit a^*$.
The $W$-invariant polynomials are described by means of these 
weights as follows.
\begin{prop} (\cite{bou}, p. 188) The $W$-invariant 
functions $q_1,\dots ,q_l$ defined by 
$q_i(H)=\sum_{\varphi\in W}e^{\sqrt{-1}\lambda_i(\varphi (H))}$
generate the set $P(\mathit a)^W$ of 
$W$-invariant trigonometric polynomials.
\end{prop}

Real valued generating set, $r_1,\dots ,r_l$, can be extracted by taking
real or imaginary parts of $q_i$. The real parts are even functions
and the imaginary parts are odd functions. $W$-invariant homogeneous
polynomials can be found by expanding each trigonometric polynomial,
$r_i$, into its Taylor series 
$r_i=\sum_{m=0}^\infty p_i^m$ about the origin, 
where the $p_i^m$ is homogeneous of 
degree $m$. One can find the independent generators of the set 
$S(\mathit a)^W$ of $W$-invariant homogeneous 
polynomials as follows \cite{da}.

By formula (8) of \cite{da}, there exist positive
integers $m_1,\dots ,m_l$ such that
\[
dp^{m_1}_1\wedge\dots\wedge dp_l^{m_l}=c
\big(\Pi_{\{\alpha\in\Delta^+\}}\alpha\big) dH,
\]
where the set $\Delta^+$ is defined by those roots which are positive
on a fixed Weyl chamber, furthermore, $dH$ is the $l$-form giving 
the Euclidean volume element on  $\mathit k$. Then one has:
\begin{prop} \cite{da} The polynomials 
$p_1^{m_1},\dots ,p_l^{m_l}$
are algebraically independent and generate $S(\mathit a)^W$. 
Furthermore,
$p^{m_i}_i(H)=\sum_{\varphi\in W}\lambda^{m_i}_i(\varphi (H))$
and $\sum m_i$ is equal to the number of positive roots.
\end{prop}
The degrees $m_i$ of the generating Chevalley polynomials 
are well known (cf. \cite{wa}, p. 144,
\cite{bou}, \cite{ch}, \cite{boch}, \cite{ko1,ko2}). Their complete 
list is as follows.

\newpage

\centerline{Table 3.1}

\centerline{Irreducible Weyl groups and their Chevalley degrees}
\centerline{
\begin{tabular}[t]{|c|c|c|c|}\hline
diagram &  rank & $W$ & $m_i$\\\hline
$A_l$ & $l$ & $\mathbf A_l$ & $2,3,4,\dots ,l+1$ \\
$B_l$ & $l$ & $\mathbf B_l$ & $2,4,\dots ,2l$ \\
$C_l$ & $l$ & $\mathbf B_l$ & $2,4,\dots ,2l$ \\
$BC_l$ & $l$ & $\mathbf B_l$ & $2,4,\dots ,2l$ \\
$D_l$ & $l$ & $\mathbf D_l$ & $2,4,\dots ,2l-2,l$ \\
$E_6$ & $6$ & $\mathbf E_6$ & $2,5,6,8,9,12$ \\
$E_7$ & $7$ & $\mathbf E_7$ & $2,6,8,10,12,14,18$ \\
$E_8$ & $8$ & $\mathbf E_8$ & $2,8,12,14,18,20,24,30$ \\
$F_4$ & $4$ & $\mathbf F_4$ & $2,6,8,12$ \\
$G_2$ & $2$ & $\mathbf G_2$ & $2,6$ \\
\hline
\end{tabular}
}
\medskip
Though the degrees $m_i$ are known,
it is not clear that which $\lambda_i$ do they belong to. From this 
point of view this description of the generators can be considered
only as a semi-explicit description. 
More explicit description can be
given by considering the individual 
Weyl groups separately. 
We do not go
into the details of such individual investigations here. 
An other imperfection is that,
up to the
knowledge of the author, no explicit description of the extended
polynomials is known in the literature so far.

\medskip

\noindent{\bf Constructing by Chevalley's
polynomials.}
Next, the above $W$-invariant homogeneous polynomials 
are used to construct irreducible Berwald metrics. 
First pick a Cartan subspace 
$\mathit a$ and  
a $\mathcal{H}_p$-invariant inner product $\langle\, ,\,\rangle$ 
on $T_p(M)$. The latter is uniquely determined up 
to a constant factor in
the considered irreducible cases. Choose also a rectangular
coordinate system $(y_1,\dots ,y_n)$ defined by an orthonormal
basis on $T_p(M)$. In what follows, the constructions are
described for absolute- and proper positive-homogeneous norms 
separately.

\noindent{\it (A) Constructing the absolute homogeneous solutions.}
In this case we choose a 
$W$-invariant homogeneous polynomial $P(H)$ of even order, $2k$, on
$\mathit a$. These 
polynomials are required to be strictly positive on the punctured space 
$\mathit a\backslash 0$. If $P(H)$ is not strictly positive yet, 
this requirement can be furnished by adding
the term $c\langle H,H\rangle^k$ to $P(H)$, where 
$c > max|P^{-}_{/S}|$ and
$P^{-}_{/S}$ means the restriction of the negative part
$P^{-}(H)=(1/2)(P(H)-|P(H)|)$ of $P(H)$ onto the unit sphere $S$
around the origin of the Euclidean space $\mathit a$. Then, 
by the Chevalley Theorem, the
polynomial $Q(H)=c\langle H,H\rangle^k+P(H)$ extends to a 
$\mathcal{H}_p$-invariant
homogeneous polynomial $\tilde{Q}(X)$ of order $2k$ 
on $T_p(M)$ such that 
$\tilde{Q}^{1/k}$ 
is a $C^\infty (\mathit a\backslash 0)$,
homogeneous function of order $2$.

As in (\ref{L1}), consider the $0$-homogeneous matrix field
${1\over 2}\partial_i\partial_j\tilde Q^{1\over k}$. 
If this field is not positive
definite and $\delta$ is the minimum of
negative eigenvalues of this field  
then the norm
\begin{equation}
L(X)=(\gamma |X|^2+
\tilde Q^{1/k}(X))^{1\over 2},
\label{Labshom}
\end{equation}
where $\gamma > |\delta |$,
is obviously a $\mathcal{H}_p$-invariant Finsler norm 
of class $C^\infty$.

This construction can be generalized by considering the norms
\begin{equation}
L(X)=(\gamma |X|^2+\sum c_i
\tilde Q_i^{1/k_i}(X))^{\frac{1}{2}},
\label{Labshom2}
\end{equation}
where $\gamma > |\delta |$ are defined with respect to the function
$\sum c_i\tilde Q_i^{1/k_i}$. 
If this latter sum consists only finite many terms, the constants
can be arbitrary positive or negative numbers. 
The $L$ is always a function
of class $C^\infty$. 

By using infinite sums, one can construct
solutions of class $C^m$ by choosing constants $c_i$ such that
$\sum c_i\tilde Q_i^{1/k_i}$ is convergent with 
respect to the $C^m$-norm
on the unit sphere $S$. This technique reveals a considerably wider
class of solutions, however, the examples introduced so far have
absolute homogeneous Finsler functions.

\noindent{\it Constructing the positive homogeneous Finsler functions.}
Now let $k=2s+1$ be an odd number and $P_{2k}$ resp. $P_k$ be 
$W$-invariant homogeneous polynomials of degree $2k$ resp. $k$. By
the discussions given below about the explicit form of $W$-invariant
polynomials, such odd numbers and polynomials exist exactly
on those symmetric spaces where the Weyl group does not contain the 
reflection $-id$ in the origin. In this case there exist $W$-invariant
homogeneous polynomials $P_k$ of odd order $k=2s+1$. The
polynomial $P_{2k}$ can
be chosen, for instance, as the product of two such polynomials. 

By adding a 
suitable correction term $c\langle H,H\rangle^k$, the polynomial
$R_{2k}(H)=c\langle H,H\rangle^k+P_{2k}(H)$ 
becomes strictly positive on 
$\mathit a\backslash 0$ and the extended function 
$\tilde{R}_{2k}^{1/2}$ 
is a positive, $C^\infty$, and absolute homogeneous 
function of degree $k$ on $T_p(M)$.

In the next step consider the positive homogeneous function
\begin{equation}
Q^*=(\tilde{R}^{1\over 2}_{2k}+\tilde{P}_k)^2=
(\tilde{R}_{2k}+\tilde{P}_k^2)+2
\tilde{R}_{2k}^{1\over 2}\tilde{P}_k
\label{L+hom1}
\end{equation}
of even degree $2k$. Note that the last term is a skew-symmetric 
function while the first two terms are absolute homogeneous. 
If the $Q^*$ is
not strictly positive on $T_p(M)\backslash 0$, then 
$\tilde{Q}(X)=d\langle X,X\rangle^k+Q^*(X)$ 
will have this property for all $d > 0$.

The construction is completed by the Finsler norm
\begin{equation}
L(X)=(\gamma |X|^2+
\tilde Q^{1/k}(X))^{1\over 2}.
\label{L+hom2}
\end{equation}
These solutions are just positive- but 
not absolute-homogeneous Finsler functions.
 
The most general solutions can be defined by formula (\ref{Labshom2}), 
where in the sum one can use both the absolute and positive
homogeneous functions constructed above.

\medskip
\noindent{\bf Classification theorems.}
These theorems can be established by combining the previous statements
with the well known Cartan-list of irreducible 
symmetric manifolds. The table below provides the complete list along
with the restricted Weyl groups belonging to these manifolds (cf.
\cite{he}, p. 532-534, \cite{kn}, p. 714, \cite{wa}, p. 30-32).

\medskip
\centerline{Table 3.2}
\centerline{
Irreducible symmetric manifolds and their restricted Weyl groups}
\centerline{
{\scriptsize
\begin{tabular}[t]{|l|l|l|l|} 
\hline
non-compact 
& compact & $W_l$ & $dim(M)$
\\
\hline
$SL(n,\mathbf R)/SO(n)$ &
$SU(n)/SO(n)$ & $\mathbf A_{n-1}$ & ${1\over 2}(n-1)(n+2)$\\
$SU^*(2n)/Sp(n)$ &
$SU(2n)/Sp(n)$ & $\mathbf A_{n-1}$ & $(n-1)(2n+1)$\\
$SU(p\geq q)/S(U_p\times U_q)$ &
$SU(p+q)/S(U_p\times U_q)$ & $\mathbf B_q$ & $2pq$\\
$SO_0(p,p)/SO(p)\times SO(p)$ &
$SO_0(p+p)/SO(p)\times SO(p)$ & $\mathbf D_p$ & $p^2$ \\
$SO_0(p>q)/SO(p)\times SO(q)$ &
$SO_0(p+q)/SO(p)\times SO(q)$ & $\mathbf B_q$ & $pq$ \\
$SO^*(2n)/U(n)$ &
$SO(2n)/U(n)$ & $\mathbf B_{[{1\over 2}n]}$ & $n(n-1)$\\
$Sp(n,\mathbf R)/U(n)$ &
$Sp(n)/U(n)$ & $\mathbf B_n$ & $n(n+1)$\\
$Sp(p\geq q)/Sp(p)\times Sp(q)$ &
$Sp(p\geq q)/Sp(p)\times Sp(q)$ & $\mathbf B_q$ & $4pq$ \\
$SL(n+1,\mathbf C)/SU(n+1)$ &
$\times_2 SU(n+1)/SU(n+1)$ & $\mathbf A_{n\geq 1}$ & $n(n+2)$\\
$SO(2n+1,\mathbf C)/SO(2n+1)$ &
$\times_2 SO(2n+1)/SO(2n+1)$ & $\mathbf B_{n\geq 2}$ & $n(2n+1)$\\
$Sp(n,\mathbf C)/Sp(n)$ &
$Sp(n)\times Sp(n)/Sp(n)$ & $\mathbf B_{n\geq 3}$ & $n(2n+1)$\\
$SO(2n,\mathbf C)/SO(2n)$ &
$SO(2n,)\times SO(2n)/SO(2n)$ & $\mathbf D_{n\geq 4}$ & $n(2n-1)$\\
$\mathcal E_{6(6)}/Sp(4)$
& $\mathcal E_{6(-78)}/Sp(4)$ & $\mathbf E_6$ & $42$\\
$\mathcal E_{6(2)}/SU(6)\times SU(2)$
& $\mathcal E_{6(-78)}/SU(6)\times SU(2)$ & $\mathbf F_4$ 
& $40$\\
$\mathcal E_{6(-14)}/SU(10)\times\mathbf R$
&$\mathcal E_{6(-78)}/SU(10)
\times\mathbf R$ & $\mathbf B_2$ & $32$\\
$\mathcal E_{6(-26)}/\mathcal F_4$
& $\mathcal E_{6(-78)}/\mathcal F_4$ & $\mathbf A_2$ & $26$ \\
$\mathcal E_{7(7)}/SU(8)$
&$\mathcal E_{7(-133)}/SU(8)$ & $\mathbf E_7$ & $70$\\
$\mathcal E_{7(-5)}/SO(12)\times SU(2)$
&$\mathcal E_{7(-133)}/SO(12)\times SU(2)$ & $\mathbf F_4$ 
& $64$\\
$\mathcal E_{7(-25)}/\mathcal E_6\times\mathbf R$
&$\mathcal E_{7(-133)}/\mathcal E_6\times\mathbf R$ & $\mathbf B_3$ & 
$54$\\
$\mathcal E_{8(8)}/SO(16)$
&$\mathcal E_{8(-240)}/SO(16)$ & $\mathbf E_8$ & $128$\\
$\mathcal E_{8(-24)}/\mathcal E_7\times SU(2)$
&$\mathcal E_{8(-248)}/\mathcal E_7\times SU(2)$ & $\mathbf F_4$  
& $112$\\
$\mathcal F_{4(4)}/Sp(3)\times SU(2)$
&$\mathcal F_{4(-52)}/Sp(3)\times SU(2)$ & $\mathbf F_4$ & 
$28$\\
$\mathcal F_{4(-20)}/SO(9)$ &
$\mathcal F_{4(-52)}/SO(9)$ & $\mathbf A_1$ & $16$\\
$\mathcal G_{2(2)}/SU(2)\times SU(2)$
&$\mathcal G_{2(-14)}/SU(2)\times SU(2)$ & $\mathbf G_2$ & $8$ \\
$\mathcal E_6^C/\mathcal E_6$
&$\mathcal E_6\times \mathcal E_6/\mathcal E_6$ & $\mathbf E_6$ &$78$\\
$\mathcal E_7^C/\mathcal E_7$
&$\mathcal E_7\times \mathcal E_7/\mathcal E_7$ & $\mathbf E_7$ 
& $133$\\
$\mathcal E_8^C/\mathcal E_8$
&$\mathcal E_8\times \mathcal E_8/\mathcal E_8$ & $\mathbf E_8$ &$248$\\
$\mathcal F_4^C/\mathcal F_4$
&$\mathcal F_4\times \mathcal F_4/\mathcal F_4$ & $\mathbf F_4$ & $52$\\
$\mathcal G_2^C/\mathcal G_2$
&$\mathcal G_2\times \mathcal G_2/\mathcal G_2$ & 
$\mathbf G_2$ & $14$\\
\hline 
\end{tabular}
}}
\medskip\medskip

In this table
formula 
$SU(p\geq q)$ denotes the group  
$SU(p,q)$ such that  
$p\geq q$ holds. The groups   
$SO_0(p\geq q)$ and  
$Sp(p\geq q)$ are similarly defined. The Cartesian product $G\times G$
is sometimes written in the form $\times_2G$. Formula 
$\mathbf D_{n\geq 4}$
indicates that in the line of the formula the relation    
$n\geq 4$ is supposed to be satisfied. The meaning of formulas
$\mathbf A_{n\geq 1}$,
$\mathbf B_{n\geq k}$ is similar.

The second table describes 
the complete list of irreducible globally symmetric
Riemannian manifolds (cf. \cite{he}, p. 516-518). Non-Riemannian Berwald
metrics can be constructed exactly on those having $rank\geq 2$.
The symmetric spaces of $rank 1$ are exactly those having
the restricted Weyl group $\mathbf A_1$ or $\mathbf B_1$. By this
principle one can sort out the well known list
\begin{eqnarray}
\label{rank1}
SU(p,1)/S(U_p\times U_1)\, ,\,
SU(p+1)/S(U_p\times U_1) \, , \\
SO_0(p,1)/SO(p)\, ,\,
SO_0(p+1)/SO(p)\, ,\nonumber\\
Sp(p,1)/Sp(p)\times Sp(1)
\, ,\, Sp(p,1)/Sp(p)\times Sp(1)\, ,\nonumber\\
\mathcal F_{4(-20)}/SO(9)\, ,\,
\mathcal F_{4(-52)}/SO(9)\nonumber
\end{eqnarray}
of simply connected and complete symmetric spaces of rank 1.

Finally, Table 3.1 reveals
that there exist skew-symmetric $W$-invariant
polynomials if and only if $W$ is either
$\mathbf A_p$ with $p\geq 2$,or, it is one of the groups 
$\mathbf D_{2k+1}$ and $\mathbf E_6$.
An irreducible Riemannian connection can be metrized by irreversible
($d(x,y)\not =d(y,x)$) Berwald metric 
if and only if it is symmetric and its Weyl group is
one of these groups. By summing up we get
\begin{theorem}{\bf (Constructive Main Theorem 1)} 
(A) An irreducible connection of Riemannian
type is Berwald metrizable by non-Riemannian Finsler metrics if and
only if it is symmetric of rank$\geq 2$. Table 3.2 lists all the 
simply connected, complete irreducible Riemannian symmetric spaces. 
By leaving out from this list
the rank 1 symmetric cases listed in (\ref{rank1}), 
one gets the complete
list of those simply connected, complete 
irreducible Riemannian manifolds whose Riemannian connection can be
metrized also by non-Riemannian Berwald metrics.

The Berwald metrics on these manifolds are explicitly constructed in
formulas (\ref{L1})-(\ref{L+hom2}).

(B)The only 
simply
connected, complete 
irreducible Riemannian manifolds whose Riemannian connection can be
metrized also by irreversible Berwald metrics are the following ones:
\begin{eqnarray}
\label{skew}
SL(n,\mathbf R)/SO(n)\, ,\,
SU(n)/SO(n) , \quad n\geq 3 ,\\
SU^*(2n)/Sp(n) \, ,\,
SU(2n)/Sp(n) , \quad n\geq 3 ,\nonumber\\
SO_0(p,p)/\times_2 SO(p) \, ,\,
SO_0(2p)/\times_2 SO(p) , \quad p=2k+1 ,\nonumber\\
SL(n+1,\mathbf C)/SU(n+1) \, ,\,
\times_2 SU(n+1)/SU(n+1) , \quad n\geq 2 ,\nonumber\\
SO(2n,\mathbf C)/SO(2n) \, ,\,
\times_2 SO(2n)/SO(2n) , \quad n=2k+1 ,\nonumber\\
\mathcal E_{6(6)}/Sp(4)
\, , \,\mathcal E_{6(-78)}/Sp(4)\,\, ,\,\,
\mathcal E_{6(-26)}/\mathcal F_4 \,\, ,\,\,
\mathcal E_{6(-78)}/\mathcal F_4 \, ,\nonumber\\
\mathcal E_6^C/\mathcal E_6\, ,\,
\mathcal E_6\times \mathcal E_6/\mathcal E_6.\nonumber 
\end{eqnarray}

Explicit proper positive homogeneous Finsler functions 
on these manifolds are described
in formulas (\ref{L+hom1}) and (\ref{L+hom2}).
The Riemannian connections on these manifolds allow also reversible
($d(x,y)=d(y,x)$) non-Riemannian Berwald metrics. 
The other manifolds mentioned in (A)
allow only reversible Berwald metrics.

\end{theorem}

\section{Reducible Berwald metrics}

For reducible holonomy groups
the above construction 
is combined with the {\it deRham decomposition} of Riemannian
connections. Let
\begin{equation}
T_p(M)=V_0\oplus S_1\oplus\dots\oplus S_l\oplus U_1\oplus\dots\oplus U_k
\label{irrdec}
\end{equation}
be the irreducible decomposition of the tangent space at $p$, with
respect to the holonomy group $\mathcal{H}_p$, 
where $V_0$ is the maximal 
point-wise fixed subspace, and the subspaces $S_i$ resp. 
$U_i$ belong to irreducible symmetric resp. to non-symmetric spaces.
The integral manifolds $M_0,M_{Si},M_{Ui}$ corresponding to the deRham
decomposition are globally defined for simply connected and complete
$M$. Otherwise they are only locally defined. By the Berger-Simons
holonomy theorem, a Berwald metric
induces a Minkowski metric
on $M_0$, the above described metrics on $M_{Si}$, 
and Riemannian metrics
on the manifolds $M_{Ui}$.

In order to construct all the Berwald metrics in the case of 
reducible Berwald
connections, the concepts of Cartan subalgebras and 
Weyl groups are introduced for arbitrary Riemannian connections. 
A {\it Cartan subspace},
$\mathit c\subset T_p(M)$, is defined by the direct sum
\begin{equation}
\mathit c=V_0\oplus \mathit c_{s1}\oplus\dots\oplus \mathit c_{sl}\oplus
\mathit c_{u1}\oplus\dots\oplus \mathit c_{uk},
\label{Cart}
\end{equation}
where $\mathit c_{si}$ is a Cartan subalgebra on 
$S_i$ and $\mathit c_{ui}$ is an
arbitrary 1-dimensio\-nal subspace in $U_i$. The Weyl group $W_0$
on $V_0$ consists only of the identity map and the Weyl group 
$W_{ui}$ acting on $\mathit c_{ui}$ contains only two elements, 
the identity map
and the central symmetry 
with respect to the origin. Then the {\it Weyl group}, 
$W$, acting on
$\mathit c$ is defined by the direct product
\begin{equation}
W=W_0\times W_{s1}\times\dots\times W_{sl}\times
W_{u1}\times\dots\times W_{uk},
\label{W}
\end{equation}
where $W_{si}$ is the usual Weyl group acting on the Cartan subalgebra 
$\mathit a_{si}=\mathit c_{si}$ 
of the corresponding symmetric Lie algebra.

This Weyl group is obviously a finite reflection 
group and both the theory
of Chevalley and the invariant theory applies to it. By summing up
we have

\begin{theorem}{\bf (Constructive Main Theorem 2)} (A) Any orbit 
$\mathcal{H}_p(X)$ 
intersects any
Cartan flat $\mathit c$ at a 
non-empty set of finite many points. The finite
Weyl group $W$ acts transitively on these intersection points.

(B)
All the $\mathcal{H}_p$-invariant Finsler norms, $L_p(X)$, defined
on $T_p(M)$ can be represented by averaging, 
$L_p(X)=\int_{\mathcal H_p} L_p^*(g(X))dg$, 
of Finsler norms $L^*_p$. Much
more explicit technique is 
offered by the above statements, however, since the norms
$L_p(X)$ in the question can be constructed by extending $W$-invariant
Finsler norms $L_{\mathit c}$, defined on a fixed  
Cartan subspace $\mathit c$, onto $T_p(M)$.

Any $W$-invariant norm $L_{\mathit c}$ can be uniquely extended to a 
continuous $\mathcal{H}_p$-invariant Banach norm 
defined on $T_p(M)$ by the formula
$L_p(X)=L_{\mathit c}(Y)$, where $Y$ is an arbitrary vector from 
$\mathit c\cap \mathcal{H}_p(X)$. There arise difficult problems 
regarding the smoothness and strong convexity 
of the extended norm, however. One can
avoid these difficulties and be able to construct appropriate Finsler
norms by means of the Chevalley polynomials, or in general, using
$C^\infty$ and $W$-invariant Finsler norms. 
Then the
extended Finsler norms appear in the form
\begin{equation}
L_p(X)=(\gamma |X|^2+\sum c_i\tilde Q_i^{1/k_i}(X))^{1\over 2},
\end{equation}
where the functions $Q_i$ can be even or 
odd, like the functions used in formulas 
(\ref{Labshom2}) and (\ref{L+hom2}).
The $C^\infty$ Finsler norms can be represented in the form
\begin{equation}
L_p(X)=(\gamma |X|^2+\tilde L_{/\mathit c}^2(X))^{1\over 2},
\end{equation}
where $L_{/\mathit c}$ is a $C^\infty$ Finsler norm defined on 
$\mathit c$ and
$\tilde{L}_{/\mathit c}$ is its extension to $T_p(M)$.

In both formulas the correction term $\gamma |X|^2$ makes up 
the norm to a strongly convex one.

The parallel extension of 
$L_p$ defines a Berwald metrization of the considered linear 
connection.

(C) The Riemannian connection of a reducible Riemannian manifold
$M=M_0\times M_1\times\dots M_k$ can be metrized by irreversible
Berwald metrics if and only if the Euclidean factor $M_0$ is non-trivial
or at least one of the other factors is equal to one of the manifolds
listed in (\ref{skew}). 

For a reducible Riemann connection 
belonging to a Riemann metric $g$,
the most simple non-Riemannian Berwald metrics are those having the
Finsler function
\begin{equation}
\label{simpL}
L(X)=\sqrt{c_1|X|^2+c_2
(|X_0|^{2p}+|X_1|^{2p}+\dots +|X_{1+s+k}|^{2p})^{1/p}},
\end{equation}
where $c_1,c_2 >0\, ,\, p >1$ and $X=X_0+\dots X_{1+s+k}$ corresponds
to the de Rham decomposition. Arbitrary positive scaling factors
can be applied also before the terms $|X_i|^{2p}$. 
\end{theorem}

Notice that the function
$(|X_0|^{2p}+|X_1|^{2p}+\dots +|X_{1+s+k}|^{2p})^{1/p}$
defines a norm for which the tensor $g^{(p)}_{ij}$ is positive semi
definite (proved by Minkow\-ski), 
therefore, the norm (\ref{simpL}) is a Finsler
norm because of the correction term $|X|^2=g(X,X)$. This latter
statement proves that any reducible Riemannian connection can be
metrized by non-Riemannian Berwald metrics.

\medskip
\noindent{\bf Remark.} The only explicit non-Riemannian 
Berwald metrics mentioned in \cite{sz1} are the metrics described 
in (\ref{simpL}) (more precisely,
they are introduced with the specific constants 
$c_1=c_2=1$). Notice that the
correction term $|X|^2$ adjusts
the convex Minkowski norm determined by 
the second term to a strongly convex one.
\medskip

One should be precautious about the definition of deRham decomposition
of Berwald manifolds. The Cartesian product concerns only the
manifold 
($M=M_0\times M_1\times\dots\times M_k$)
and the Berwald connection
($\nabla=\nabla_0\times \nabla_1\times\dots\times \nabla_k$), the 
metric,
however, may not be the Cartesian product of the metrics defined on the 
factor manifolds. It is defined by the weaker requirement demanding 
that the Finsler function should be invariant
with respect to the holonomy group. 

This requirement allows infinitely
many inequivalent solutions for defining Finsler functions on the
product manifold $M$ by means of the Finsler functions $L_i$ 
defined on the
factor manifolds $M_i$. The Cartesian product is just one of the 
options. In order to make a clear distinction between Cartesian
product and the product applied in this paper, the
latter is called {\it holonomy-invariant, or, 
perturbed Cartesian product of the metrics}. 
Then the deRham decomposition of Berwald spaces can be stated as 
follows.

\begin{theorem}\label{ST}
{\bf (Structure Theorem)}
Let $(M,\nabla ,L)$ be a connected, simply connected, 
and complete Berwald
manifold. Then the affine manifold $(M,\nabla )$ decomposes into
the Cartesian product
$(M=M_0\times M_1\times\dots\times M_k ,
\nabla=\nabla_0\times \nabla_1\times\dots\times \nabla_k)$,
where 
$(M_0,\nabla_0,L_0)$
is the maximal Minkowskian factor and the irreducible factors 
$(M_i,\nabla_i,L_i)$
are either Riemannian manifolds or non-Riemannian affine symmetric
Berwald manifolds. The metric $L$ is a holonomy-invariant product
of the metrics $L_j$ defined on the factor manifolds.
\end{theorem}

Finally, we describe a straightforward generalization 
of Hano's Theorem concerning the decomposition of group of isometries,
induced by the de Rham decomposition.
\begin{theorem}
{\bf (Generalized Hano Theorem \cite{sz1})} Let 
$M=M_0\times\dots\times M_{1+s+k}$
be the de Rham decomposition of a complete simply connected 
Berwald manifold. If $U_0(M_i)$ and $J_0(M_i)$ denote the unit 
component of the affine and the isometry group on $M_i$ respectively,
then $U_0(M)=U_0(M_0)\times\dots\times U_0(M_{1+s+k})$ resp.
$J_0(M)=J_0(M_0)\times\dots\times J_0(M_{1+s+k})$.
\end{theorem}

\section{Cartan-symmetric Finsler manifolds}

A diffeomorphism, $\varphi:M\to M$, is an isometry on a
Finsler manifold $(M,L)$ if it preserves the Finsler function.
By the classical Dantzig-van der Waereden Theorem 
(cf. \cite{kn}, vol. I, chapter I, Theorem 4.7) asserting:
{\it ``The group $G$ of isometries of a 
connected, locally compact
metric space $M$ is locally compact with respect 
to the compact-open topology''},
and the Montgomery-Zippin Theorem 
(cf. \cite{kn}, vol. I, chapter I, Theorem 4.6) asserting:
{\it ``A locally compact effective transformation group $G$ of
$C^1$-diffeomorphisms acting on a connected manifold $M$ of class $C^k$
is a Lie group and the mapping $G\times M\to M$ is of class $C^k$''},
the group of isometries on a connected Finsler manifold form a
Lie group. Strictly speaking, these theorems prove the statement for
absolute homogeneous Finsler functions. For positive 
homogeneous Finsler functions consider the
metric, $d^*$, defined by the function $L^*(X)=L(X)+L(-X)$. Then 
the $G$ is a closed subgroup of $G^*$ defined for $d^*$. Thus both 
groups are Lie groups.

A Finsler manifold is called
{\it Cartan-symmetric} if there exist involutive isometry 
$\sigma_p$, 
for each point $p\in M$,
for which the $p$ is an isolated fixpoint. 
It is obvious that
such $\sigma_p$ induces $-id$ on the tangent space $T_p(M)$,
therefore, Cartan-symmetric Finsler manifolds have reversible
metrics and geodesics. 

The geodesic symmetries 
$s_p: exp(y)\to exp(-y)$
on Berwald manifolds having symmetric Riemannian connections 
are investigated in \cite{dh1, dh2}. 
They are isometries (Cartan-symmetries) 
if and only if the symmetric B-manifold has
absolute homogeneous Finsler function. A general
Berwald manifold with symmetric Riemannian connection is called
{\it affine symmetric Berwald manifold}, since such a manifold with 
proper positive homogeneous Finsler function is not 
Cartan-symmetric. Note that the Cartan symmetric manifolds are defined
among the most general Finsler manifolds and not just among 
the Berwald manifolds. Next we prove, however, that these
general Cartan-symmetric 
Finsler manifolds are Berwald manifolds, having affine symmetric 
connection $\nabla$ and an absolute homogeneous Finsler function.

In the following considerations we suppose that 
the symmetries $\sigma_p$ are globally
defined. By passing to the universal covering space,
one can also suppose that the $M$ is simply connected. The above
statement is established via several steps.

\noindent{\it (A) Cartan-symmetric manifolds are homogeneous.} 
Indeed, for two points
$p$ and $q$ which can be connected by a geodesic segment $s(p,q)$, the 
central symmetry $\sigma_m$, where $m$ is the midpoint of $s(p,q)$, 
corresponds $p$ and $q$ to each other. Note that this argument 
exploits that Cartan symmetric manifolds have reversible geodesics. 
Considering two arbitrary
points $p$ and $q$, connect them by a sequence
of geodesics. Then the composition of the central symmetries 
defined by the 
midpoints of these geodesics sends the starting point $p$ to the 
endpoint $q$. This proves the statement completely. 

\noindent{\it (B) The group of isometries is a Lie group} by the 
argument described at the very beginning of this section. 
{\it Furthermore, the isotropy group  $H_p$ of isometries 
fixing a point $p$ is compact.} 
(This latter statement can be established by observing
that the indicatrix in the tangent space $T_p(M)$ is invariant
under the actions of the induced maps $\varphi_*$, 
for all $\varphi\in H_p$.)
It follows that globally Cartan-symmetric manifolds are complete ones.

\noindent{\it (C) The homogeneous space $J(M)/H_p$ is symmetric.} 
In fact, let $G_\sigma\subset J$ be the subgroup of isometries 
generated by  
both algebraic and topological closure of the set of 
Cartan-symmetries. Because of closedness, it is a
Lie group, moreover, it is a normal subgroup of $J$. Let 
$G_{\sigma 0}$ be the subgroup generated by those isometries
which can be expressed as products of even number of Cartan-symmetries.
This subgroup is in the unit component of $G_\sigma$.
For let $c_1(t),\dots ,c_{2k}(t)$ be continuous curves connecting
the point $p=c_1(0)=\dots =c_{2k}(0)$ with the points $p_i=c_i(1)$.
Then the curve 
$S(t)=\sigma_{c_ 1(t)}\circ\dots\circ\sigma_{c_{2k}(t)}$
in $G_{\sigma 0}$ connects
$S(1)=\sigma_{p_ 1}\circ\dots\circ\sigma_{p_{2k}}$
with $id=S(0)$. For any two points $p$ and $q$ can be connected
by even number of geodesics, 
the group $G_{\sigma 0}$ is a transitive and normal subgroup of
the isometries.

Let $H_{\sigma p0}\subset G_{\sigma 0}$
be the subgroup whose elements fix a point $p\in M$. Then $M$ is a 
homogeneous space of the form 
$M=G_\sigma /H_{\sigma p}=G_{\sigma 0}/H_{\sigma p0}$.
One can construct all of the $G_\sigma$-invariant Cartan-symmetric
Finsler metrics on $M$ by extending the 
$H_{\sigma p}$-invariant Finsler norms,
defined on $T_p(M)$, by the $G_\sigma$ onto $M$. 
Since the $H_{\sigma p}$ is compact, there exists also a
$H_{\sigma p}$-invariant inner product on $T_p(M)$ which extends into
a $G_\sigma$-invariant Riemannian metric of $M$. 
(See more about these constructions in \cite{sz3} and/or 
Theorems 1.3 and 1.4 in \cite{dh1}.)
 
Since the isometry 
group of this Riemannian metric contains the 
central symmetries $\sigma_q$ for all $q\in M$,
it defines a symmetric Riemannian space \cite{kn}. (The standard proof
of this statement is as follows. 
Let $G$ be the unit component of isometries on
this Riemannian manifold and $H_p\subset G$ 
be the isotropy subgroup fixing
the point $p$. Then $M=G/H_p$. Consider the involutive authomorphism
$\sigma (g)=\sigma_p\circ g\circ\sigma_p$ on $G$. Then 
$K_{\sigma 0}\subset H_p\subset K_\sigma$, where $K_\sigma$ is the set
of fixpoints of the $\sigma$ and $K_{\sigma 0}$ is its unit component.
Therefore, the $(G,H_p)$ is a Riemannian symmetric pair, proving the
statement completely.)

On irreducible affine symmetric Berwald
manifolds the  
unit component, $J_0$, of the isometries 
is a simple group (cf. Cartan's theorem), in which the
$G_{\sigma 0}\subset J_0$ is a normal subgroup. 
Thus, $G_{\sigma 0}=J_0,\, H_{\sigma p0}=H_{p0}$ hold. 
Furthermore, there is a natural
identification between the holonomy group at a point $p$ and the 
isotropy group $H_{p0}$.
These arguments and the de Rham decomposition described earlier 
establish also the Generalized Hanno Theorem completelly.

Consider the Riemannian
connection $\nabla$ of this Riemannian symmetric manifold. If
$\tau_c :T_p(M)\to T_q(M)$ 
is the parallel transport along a curve $c$ joining
$p$ and $q$ then there exist $\varphi\in G_{\sigma 0}$ satisfying
$\varphi (p)=q$ such that the induced map
$\varphi_* :T_p(M)\to T_q(M)$ is nothing but $\tau$ (cf. Theorem 3.2
of chapter XI in \cite{kn}, vol. II). Therefore
$\nabla_k L=0$ holds for any $G_\sigma$-invariant Finsler metric.
This proves that a Cartan-symmetric Finsler metric is always
an affine symmetric Berwald metric with an absolute 
homogeneous Finsler function. Thus also statement (C) is established. 

\noindent{\it (D) The converse statement, asserting that any affine 
symmetric Berwald metric having reversible Finsler function is 
Cartan-symmetric, is also true.} 
In this case the Finsler function is defined
by parallel transports of a Finsler function $L_p$, defined for a fixed
point $p$ of $M$, to the other points of the manifold. This 
construction insures, by the theorem quoted from \cite{kn} vol. II 
above, that the maps from $G_{\sigma 0}$ act as
isometries on the manifold. In order to establish this property
for any Cartan symmetry $\sigma_q$, we have to prove yet that for any
two points $p_1$ and $p_2=\sigma_q(p_1)$ the induced map 
$\sigma_{q*}$ transports $L_{p_1}$ to $L_{p_2}$. 

First note that  
$G_\sigma =G_{\sigma 0}\cup \sigma_pG_{\sigma 0}$ and 
$H_{p} =H_{p0}\cup \sigma_pH_{p0}$ 
hold, which statements follow from 
$
\sigma_{p_1}\circ\dots\circ\sigma_{2k+1}
=\sigma_p\circ\sigma_p\circ\sigma_{p_1}\circ\dots\circ\sigma_{2k+1}
$ immediately. Furthermore, there exist maps 
$\alpha_1, \alpha_2\in G_{\sigma 0}$
such that $\alpha_1(p)=p_1, \alpha_2(p)=p_2$ hold (cf. (C)), therefore,
the map $\varphi =\alpha^{-1}_2\circ\sigma_q\circ\alpha_1$ fixes $p$. 
Moreover, the $\alpha_{i*}$ transports $L_p$ to $L_{p_i}$, thus,
$\sigma_{q*}$ transports $L_{p_1}$ to $L_{p_2}$ if and only if
$L_p$ is invariant under the action of 
$\varphi_{p*}: T_p(M)\to T_p(M)$. 
Reversible Berwald-Finsler functions, $L_p$,
are invariant under the action of $\sigma_{p*}$ and $h_{p0*}$, for all
$h_{p0}\in H_{p0}$, which statement proves this invariance and the 
converse statement (D) completely.

By summing up (A) through (D), we have

\begin{theorem}
The Cartan-symmetric Finsler metrics on simply connected and complete
manifolds are exactly the a\-ffine symmetric Berwald metrics 
having absolute homogeneous (reversible) Finsler functions.
\end{theorem}

By this theorem, the classification of Berwald metrics gives a complete 
list also for Cartan-symmetric Finsler metrics.

\section{Cartan flats on Berwald manifolds}

Let $T_{\mathit c}$ be the surface described by 
the geodesics tangent to the
Cartan subspace $\mathit c=V_0\oplus_i\mathit c_i\subset T_p(M)$ at 
a $p\in M$ of the Berwald manifold $M$. 
Considering it as an affine submanifold with 
the induced affine connection $\tilde\nabla$, the $T_{\mathit c}$ 
is a Cartesian product of 
the totally geodesic flats $M_0=T_{V_0}$ and $T_{\mathit c_i}$. 
(Note that
on an irreducible non-symmetric manifold the surface $T_{\mathit c_i}$ 
is one-dimensional defined by a geodesic.) 
Therefore, the $T_{\mathit c}$ is a totally geodesic flat surface on 
which the induced metric is a Minkowski metric. 

These flats are called {\it Cartan flats}. 
The uniquely determined dimension of the flats is
equal to the rank of the Berwald manifold.

Since the parallel transports 
$\tau :T_p(M)\to T_q(M)$ 
along the curves joining
$p$ and $q$ corresponds Cartan subspaces to each other such that they 
keep also the Finsler function $L_{/\mathit c}$, any two Cartan 
flats with the induced Minkowski metrics are locally isometric.
Thus we have
\begin{theorem}
Any tangent vector $v$ is contained in at least one of the
Cartan flats. Any two Cartan flats of a Berwald manifold are locally 
isometric. By the above construction technique, 
the Minkowskian
metric on a single Cartan flat uniquely determines
the metric on the
whole Berwald manifold. 
\end{theorem}

We conclude the paper by a characterization of affine symmetric Berwald
manifolds, which reminds the theorem asserting that 
rank-one symmetric
spaces are exactly the 2-point homogeneous Riemannian metrics. 
(The precise formulation of this
latter theorem is as follows: 
An irreducible Riemannian
manifold of $dim > 1$ is rank-one symmetric if and only if the 
isometries from the unit component $J_0$ act transitively on the set
of pairs 
$(p,\mathit t)$, 
where $\mathit t$ is a geodesic and 
$(p\in\mathit t)$ is a point on it.) 
The corresponding characterization is:
  
\begin{theorem}\label{kflat}
A Berwald metric is affine symmetric if and only
if the isometries $J_0(M)$ act transitively on 
the set of pairs $(p,T_{\mathit c})$,
where $T_{\mathit c}$ is a totally geodesic Cartan flat with a
distinguished point $p\in T_{\mathit c}$.
\end{theorem}

The first unified proof, working both in the compact and non-compact
cases, of symmetry of 2-point homogeneous spaces is given in \cite{sz2}.
The proof is elementary, using only simple topological ideas. The
following ``point version'' of the theorem is established there: 
{\it If the isotropy group $J_p$ of isometries fixing the point $p$ are 
transitive on the unit sphere in the tangent space $T_p(M)$, then
$(\nabla R)_{/p}=0$ at $p$.}

Theorem \ref{kflat} immediately follows from the 
Generalized Hano Theorem and from this
theorem. In fact, the transitivity in the theorem implies that 
also the rank-one-factors in the deRham decomposition
are symmetric spaces of rank one. The irreducible factors of higher
rank are automatically symmetric spaces.

The above type of characterization of symmetric spaces
appeared first in \cite{hptt}, where the compact 
symmetric spaces of rank $k$
are identified with the so called $k$-flat homogeneous manifolds. We
should emphasize, however, that the $k$-flats are introduced 
in that paper by much
weaker assumptions, namely, by the connected closed 
k-dimensional totally geodesic flat submanifolds of Riemannian
manifolds. Note that this concept does not refer to Cartan subalgebras.
The $k$-flat homogeneity
is introduced on those manifolds where any tangent vector
$v$ is in at least one $k$-flat $\sum$ and is defined by the 
property that the isometries are transitive on the set of 
pairs $(p,\sum )$, where 
$p$ is a point on the $k$-flat $\sum$. The proof
under such weak conditions
is much more challenging and difficult. 
In \cite{hptt} the characterization is established only on compact
manifolds, by using even the classification of symmetric spaces. Still
on compact manifolds but using weaker assumptions, the 
above characterization
is established also in \cite{eo}. Proofs established both in the
compact and non-compact cases are due to Samiou \cite{sa}.

The concept of rank introduced in this paper has an apparent relation
also to the rank-concept intensely 
investigated in \cite{bbe,bbs,b,bs}. In
these papers the rank of a geodesic in a Riemannian manifold is defined
by the dimension of the vector space spanned by the parallel 
Jacobi fields along the geodesic. The rank of the Riemannian manifold
is defined by the minimum of the ranks of the geodesic. This is even 
a much more subtle consideration of the rank-concept. 
There is proved in these papers that any locally irreducible,
compact Riemannian manifold with non-positive sectional curvature
and higher rank is a locally symmetric space. 
\medskip

\noindent{\it Zolt\'an Imre Szab\'o, Lehman College of CUNY, Bronx, N.Y., 10468
and R\'enyi Institute, Budapest, Hungary.}

\noindent{\it email: zoltan.szabo@lehman.cuny.edu}

\end{document}